\documentclass[12pt]{article}
\usepackage{amsmath}
\usepackage{amssymb}
\usepackage{theorem}

\usepackage{color}

\newcommand{\eop}{\bigstar}

\newcommand{\card}[1]{{\vert #1 \vert} }

\newcommand{\llg}{{l\rm g}}

\newcommand{\cf}{{\rm cf}}

\newcommand{\initial}{\vartriangleleft}

\newenvironment{Proof}{\noindent{\bf Proof.}}{\par\bigskip} 

{\theorembodyfont{\rmfamily}\newtheorem{THEOREM}{Theorem}[section]}

{\theorembodyfont{\rmfamily}\newtheorem{Conclusion}[THEOREM]{Conclusion}}

{\theorembodyfont{\rmfamily}\newtheorem{LEMMA}[THEOREM]{Lemma}}

{\theorembodyfont{\rmfamily}\newtheorem{Main Theorem}[THEOREM]{Main Theorem}}
\newenvironment{main Theorem}{\begin{Main Theorem}} 
{\end{Main Theorem}}

{\theorembodyfont{\rmfamily}\newtheorem{Question}[THEOREM]{Question}}

{\theorembodyfont{\rmfamily}\newtheorem{Preservation Lemma}[THEOREM]{Preservation Lemma}}
\newenvironment{preservation lemma}{\begin{Preservation Lemma}}
{\end{Preservation Lemma}}

{\theorembodyfont{\rmfamily}\newtheorem{Construction Lemma}[THEOREM]{Construction Lemma}
\newenvironment{construction lemma}{\begin{Construction Lemma}}
{\end{Construction Lemma}}

{\theorembodyfont{\rmfamily}\newtheorem{Theorem}[THEOREM]{Theorem}}

{\theorembodyfont{\rmfamily}\newtheorem{Definition}[THEOREM]{Definition}}

{\theorembodyfont{\rmfamily}\newtheorem{Conventions}[THEOREM]{Conventions}}

{\theorembodyfont{\rmfamily}\newtheorem{Main Definition}[THEOREM]{Main Definition}}
\newenvironment{main definition}{\begin{Main Definition}}
{\end{Main Definition}}

{\theorembodyfont{\rmfamily}\newtheorem{Lemma}[THEOREM]{Lemma}}

{\theorembodyfont{\rmfamily}\newtheorem{Notation}[THEOREM]{Notation}}

{\theorembodyfont{\rmfamily}\newtheorem{Convention}[THEOREM]{Convention}}

{\theorembodyfont{\rmfamily}\newtheorem{Note}[THEOREM]{Note}}

{\theorembodyfont{\rmfamily}\newtheorem{Observation}[THEOREM]{Observation}}

{\theorembodyfont{\rmfamily}\newtheorem{Remark}[THEOREM]{Remark}}

{\theorembodyfont{\rmfamily}\newtheorem{Main Fact}[THEOREM]{Main Fact}}
\newenvironment{main Fact}{\begin{Main Fact}}{\end{Main Fact}}

{\theorembodyfont{\rmfamily}\newtheorem{Fact}[THEOREM]{Fact}}

{\theorembodyfont{\rmfamily}\newtheorem{Subfact}[THEOREM]{Subfact}}

{\theorembodyfont{\rmfamily}\newtheorem{Claim}[THEOREM]{Claim}}

{\theorembodyfont{\rmfamily}\newtheorem{Main Claim}[THEOREM]{Main Claim}}
\newenvironment{main claim}{\begin{Main Claim}}{\end{Main Claim}}

{\theorembodyfont{\rmfamily}\newtheorem{Corrolary}[THEOREM]{Corrolary}}

{\theorembodyfont{\rmfamily}\newtheorem{Subclaim}[THEOREM]{Subclaim}}

{\theorembodyfont{\rmfamily}\newtheorem{Corollary}[THEOREM]{Corollary}}

{\theorembodyfont{\rmfamily}\newtheorem{Example}[THEOREM]{Example}}

{\theorembodyfont{\rmfamily}\newtheorem{Proposition}[THEOREM]{Proposition}}

{\theorembodyfont{\rmfamily}\newtheorem{Discussion}[THEOREM]{Discussion}}

\newenvironment{Proof of the Subfact}
{\noindent{\bf Proof of the Subfact.}}{\par\bigskip}

\newenvironment{Proof of the Theorem}
{\noindent{\bf Proof of the Theorem.}}{\par\bigskip}

\newenvironment{Proof of the Conclusion}
{\noindent{\bf Proof of the Conclusion.}}{\par\bigskip}

\newenvironment{Proof of the Observation}
{\noindent{\bf Proof of the Observation.}}{\par\bigskip}

\newenvironment{Proof of the Fact}
{\noindent{\bf Proof of the Fact.}}{\par\bigskip}

\newenvironment{Proof of the Lemma}
{\noindent{\bf Proof of the Lemma.}}{\par\bigskip}

\newenvironment{Proof of the Claim}
{\noindent{\bf Proof of the Claim.}}{\par\bigskip}

\newenvironment{Proof of the Subclaim}
{\noindent{\bf Proof of the Subclaim.}}{\par\medskip}

\newenvironment{Proof of the Main Claim}
{\noindent{\bf Proof of the Main Claim.}}{\par\bigskip}

\newcommand{\elementary}{\prec}

\newcommand{\Bbf}{\mathbb}

\newcommand{\into}{\rightarrow}

\newcommand{\rest}{\upharpoonright}  % restriction

\newcommand{\deq}{\buildrel{\rm def}\over =}
\newcommand{\BB}{{\cal B}}
\newcommand{\CC}{{\cal C}}

\newcommand{\FF}{{\cal F}}
\newcommand{\GG}{{\cal G}}

\newcommand{\KK}{{\cal K}}

\newcommand{\PP}{{\cal P}}

\title{Club guessing and the universal
models}
\author{Mirna D\v zamonja\\
School of Mathematics\\
University of East Anglia\\
Norwich, NR4 7TJ, UK
\\
\scriptsize{M.Dzamonja@uea.ac.uk}\\
\scriptsize{http://www.mth.uea.ac.uk/people/md.html}}

\date{published in {\em Notre Dame Journal of
Formal Logic}, vol. 46, no. 3, (2005), 283-300}
\begin{document}

\baselineskip=16pt
\binoppenalty=10000
\relpenalty=10000
\raggedbottom

\maketitle

\begin{abstract} We survey the use of club guessing and other pcf
constructs in the context of showing that a given partially 
ordered class of objects does not have a largest, or a universal
element. This version contains a correction to Definition 1.1.
\footnote{
The author thanks EPSRC for their support on as Advanced Fellowship.
AMS 2000 Classification: 03C55, 03E04, 03C45.

Keywords: universal models, club guessing.}
\end{abstract} 

\section{Introduction}\label{intro} A natural problem in mathematics is
the following: given some partially ordered or a quasi-ordered
set or a class, is there
the largest element in it. An aspect of this question appears
in the theory of order where one concentrates on the properties of the
set and the partial order, ignoring the properties of the individual elements of
the set. A very different view of this question is obtained when one takes
the point that it is the structure of individual elements that is of interest.
An instance of this is the question of universality. 
Here we are given a class or a set
of objects and a notion of embedding between them, and we ask if there is
an object in the class that embeds all the others. To simplify our exposition
here we shall always assume that we are working with a set of objects and we shall
discuss the smallest cardinality of the subset of that set that
has the property of embedding all the other objects in the set. We shall refer
to this question as {\em the universality problem}, the class and the embedding
to which the problem refers will always be clear from the context. The number
mentioned above will be then referred to as {\em the universality number}.

Instances of the universality problem have been of a continuous interest to mathematicians,
especially those studying the mathematics of the infinite- even Cantor's work on the
uniqueness of the rational numbers as the countable dense linear
order with no endpoints is a result of this type. For some more recent examples
see \cite{AB}, \cite{Komjath} and \cite{Todor}. Apart from its intrinsic interest
the universality problem has an application in
model theory, more
specifically in classification theory. There it is used to distinguish between
various kinds of unstable theories. For more on this programme the reader may
consult the introduction to \cite{DjSh 710} and some of the
main results will be mentioned in \S\ref{models}. A good source about the classical
results about universality is the introduction to \cite{KjSh 409}.
The study of the universality problem
can naturally be divided into `the positive' and `the negative' part.
On the positive side one tries to show that the universality number
has at most the given value, for example 1. Proofs here are often
explicit constructions of universal objects, or forcing constructions,
see \cite{ChKe},
\cite{Sh 100} or \cite{DjSh 614}. On the negative side one
does the opposite, showing that the universality number is at least
a given value. In this paper we shall concentrate on the negative side
of the universality problem, in particular on the instances of it that are
obtained using a specific method that has appeared as a consequence of the
discovery of the pcf theory, namely the club guessing method. 
There are several other methods
that appear in the study of the negative side of the universality problem,
notably the $\sigma$-functor of Stevo Todor\v cevi\'c \cite{Todor}, but here we shall only
concentrate on the club guessing method. This method is
due to Menachem Kojman and Saharon Shelah \cite{KjSh 409}. 
We shall start by recalling the basic principles behind it and then
give some applications including the original one from \cite{KjSh 409}
to linear orders. Let us also note that in this subject it is often
not difficult to construct a specific universe of set theory in which the desired
negative universality result holds.
For example, one can find in \cite{KjSh 409} a proof
of the fact that when one forces over a model of $GCH$ to add $\lambda^{++}$
Cohen subsets to a regular cardinal $\lambda$, then in the
resulting universe there is no universal graph
of size $\lambda$, no universal linear order of size $\lambda$, or a
universal model of any first order countable theory unstable in $\lambda$.
The point of the negative universality results obtained by the club guessing
method is that they are implications between a certain pcf statement
and the desired negative universality result, so they hold in more than just one
specifically constructed universe.

This paper is a survey of some of the existing techniques and results in this subject.
Due to the extensive literature it has not been possible to mention all
relevant results, so our apologies go to the authors of the many deserving papers
which we failed to mention.

\section{Invariants and linear orders}\label{orders}
Let us start with an easy example of a negative universality
result: we remind the reader of why it is that for an infinite cardinal
$\lambda$ there is no well ordering of size $\lambda$ to which
there is an order-preserving embedding from any well order of size $\lambda$.
The reason for this is that any well order of size $\lambda$ is ordered in a
type $\zeta<\lambda^+$ and hence cannot embed ordinals larger than $\zeta$.
This simple proof has three important elements: invariants, construction and
preservation. Specifically, to each well order we have associated an {\em invariant},
namely its order type, then we observed that the invariant is {\em preserved},
in the sense that it can only increase under embedding, and finally we
have {\em constructed} a family of well orders of size $\lambda$
where many different values of the invariant
are present (namely the ordinals in $[\lambda,\lambda^+)$), so showing that
no single well order of size $\lambda$ can embed them all.

These same principles are present in many settings, for example \cite{AB},
and in particular in the club guessing proofs. The matters of course tend to
be more complex. We shall now show the original Kojman-Shelah example of
the use of club guessing for a universality result about linear
orders, cf. \cite{KjSh 409}.

Let $\lambda$ be a regular 
cardinal and let $\KK$ be the class of all linear
orders whose size is $\lambda$. By identifying
the elements of $\KK$ that are isomorphic to each other we obtain a set of
size at most $2^\lambda$, which we shall call $\KK$ again.
We may without loss of generality assume that the universe of
each element of $\KK$ is the $\lambda$ itself. We are interested
in the universality number of $\KK$, where for our notion of the embedding
we take an injective order preserving function. For future purposes let us also
fix a ladder system $\bar{C}=\langle C_\delta:\,\delta\in S\subseteq
\lambda\rangle$, such that each $C_\delta$ is a club of the corresponding $\delta$
and $S$ is some stationary subset of the set of limit ordinals below $\lambda$.
Let $C_\delta=\langle \alpha^\delta:\,i<i^\ast(\delta)\rangle$ be the increasing
enumeration.

Every member $L$ of $\KK$ can be easily represented as a
continuous increasing union
$L=\bigcup_{j<\lambda} L_j$ of linear orders of size $<\lambda$, and such
a representation is of course not unique. Any such sequence $\bar{L}=\langle L_j:\,
j<\lambda\rangle$ is called a {\em filtration} of $L$. Next we shall define
invariants for elements of $L$, but the definition of an invariant will
depend both on the filtration $\bar{L}$, the specified ladder system
$\bar{C}$ \textcolor{red}{and an ordinal $\delta$}.

\begin{Definition}\label{invariants} Suppose that $L,\bar{L},\bar{C}$ are as above
and $\delta\in S$ is such that the universe of $L_\delta$ is $\delta$. \textcolor{red}{For $\zeta\ge\delta$}, we define
the {\em invariant} ${\rm inv}_{\bar{L},\bar{C},\textcolor{red}{\delta} }(\textcolor{red}{\zeta})$ as the set
\[
\{i<i^\ast(\delta):\,(\exists \beta
\in L_{\alpha^\delta_{i+1}}\setminus L_{\alpha^\delta_{i}})
[\{x\in L_{\alpha^\delta_i}:\,x<_L\beta\}=\{x\in L_{\alpha^\delta_i}:\,x<_L \textcolor{red}{\zeta}\}]\}.
\]
\end{Definition}

So with the notation of Definition \ref{invariants}, the invariant of
$\textcolor{red}{\zeta}$ is a subset of $i^\ast_\delta$ that codes the `reflections'
of $\textcolor{red}{\zeta}$ along the places in the filtration $\bar{\lambda}$ that are
determined by $\bar{C}$. It is easy to check that the set of $\delta$
for which the universe of $L_\delta$ is $\delta$ is a club of $\lambda$,
so since $S$ is stationary there is a club of $\delta$ for which ${\rm inv}_{\bar{L},\bar{C},\textcolor{red}{\delta} }(\textcolor{red}{\zeta})$
is well defined \textcolor{red}{for $\zeta\ge\delta$}. It is also easy to see
that for any two filtrations $\bar{L}$ and $\bar{L}'=\langle L'_j:\,j<\lambda\rangle$
of $L$ there is a club of $\delta$ such that $L_\delta=L'_\delta$, hence
the dependence of the invariant on the filtration is only up to a club.
This is not the case with its dependence on the ladder system, and in fact
only certain ladder systems are of interest to us:

\begin{Definition}\label{guessing} A ladder system
$\bar{C}=\langle C_\delta:\,\delta\in S
\rangle$ is said to be a {\em club guessing sequence} iff for every club $E$
of $\lambda$ there is $\delta\in S$ such that $C_\delta\subseteq E$.
\end{Definition}

Club guessing sequences were introduced by Shelah in \cite{Sh -g} as a tool towards
the development of pcf theory and have since found many applications in various
contexts. 
There is also a number of variants of this concept, for example a number of
interesting results about various kinds of club sequences appears in Tetsuya Ishiu's 
thesis \cite{Tetsu}. However, for the moment
we shall concentrate on the simple club guessing mentioned above. An interesting question
is when such sequences exist, and one of the most important theorems in this
vein is the following:

\begin{Theorem}\label{existenceofcg}[Shelah, \cite{Sh -g}] Suppose that $\kappa$
and $\lambda$ are regular cardinals such that $\kappa^+<\lambda$. Then there
is a club guessing sequence of the form $\bar{C}=\langle C_\delta:\,\delta\in
S^\lambda_\kappa\rangle$.
\end{Theorem}

We have used the notation $S^\lambda_\kappa$ to denote the set of $\alpha<\lambda$
whose cofinality is $\kappa$. A club guessing sequence of the form
appearing in Theorem \ref{existenceofcg} will be referred to as
an {\em $S_\kappa^\lambda$-club guessing sequence}. Note that we may
without loss of generality, by intersecting with a club of order type $\kappa$ if
necessary, assume that each $C_\delta$ has order type $\kappa$. 
In fact, a stronger theorem then \ref{existenceofcg}
was proved by Shelah, showing that in addition to this restriction
on the order type the club guessing sequence above may be
assumed to have
a square-like property, in the terminology of $\cite{DjSh 710}$ to
be a {\em truly tight $(\kappa,\lambda)$ club
guessing sequence}:

\begin{Theorem}\label{existenceofcg+}[Shelah, \cite{Sh 420}] Suppose that $\kappa$
and $\lambda$ are regular cardinals such that $\kappa^+<\lambda$. 
Then there is a stationary set $S\subseteq S^\lambda_\kappa$, 
a club guessing sequence of the form $\bar{C}=\langle C_\delta:\,\delta\in S
\rangle$ and a sequence $\bar{P}=\langle P_\alpha:\,\alpha<\lambda\rangle$
such that each $P_\alpha\subseteq \PP(\alpha)$ has size $<\lambda$ and the
sequences $\bar{C}$ and $\bar{P}$ interact in the sense that
\[
\alpha\in C_\delta\setminus \sup(C_\delta\cap\alpha) \implies C_\delta\cap\alpha
\in \bigcup_{\beta<\alpha}P_\beta.
\]
\end{Theorem}

A nicely written proof of Theorem \ref{existenceofcg}
appears in \cite{KjSh 409}. Theorem \ref{existenceofcg+} can be read off from
the conjunction of the claims in \S1 of \cite{Sh 420}, and a direct proof is
given in [\cite{Sh E12}, 1.3.(a)].
Theorem 
\ref{existenceofcg} refers to cardinals $\kappa, \lambda$ which have at least one
successor gap, $\kappa^+<\lambda$. It is natural to ask what happens at $\lambda=
\aleph_1$. Clearly, in the presence of $\clubsuit$ there is a 
$S_{\aleph_0}^{\aleph_1}$-guessing sequence. A theorem of Shelah
(Theorem III 7.1 in \cite{Sh -f}) shows that it is consistent to have $2^{\aleph_0}=
\aleph_2$ and that there is no $S_{\aleph_0}^{\aleph_1}$-club guessing sequence
(in fact the theorem shows something stronger than what is being quoted here).
This also follows from the conjunction of 
a result of Shelah in \cite{Sh 100} which shows it consistent
to have the negation of $CH$ and the existence of a
universal linear order of size $\aleph_1$,
and the Kojman-Shelah's theorem \ref{KjSh} which implies that
in a model in which $CH$ fails and there is a $S_{\aleph_0}^{\aleph_1}$-guessing sequence
there cannot be a universal linear order of size $\aleph_1$. In the other direction,
a consistency result showing how to add club guessing sequences with strong
guessing properties was introduced by
P\'eter Komj\'ath and Matthew Foreman in \cite{FoKo}, who give a cardinal-preserving
forcing which for any regular $\lambda$ and a
given stationary set $S$ in $\lambda$, keeps $S$ stationary and adds a 
a sequence $\langle C_\delta:\,\delta\in S\rangle$ such that
for every club $E$ of $\lambda$ there is a club $C$ such that for all $\alpha\in
C\cap S$ a non-empty final segment of $C_\delta$ is included in $E$.

We shall now state and give a sketch of the proof of the Kojman-Shelah theorem on
linear orders, using the notions of an invariant and a club guessing sequence
introduced above. For the case of $\kappa=\aleph_0$ we only need to use
a $S^\lambda_\kappa$ club guessing sequence, while the case $\kappa>\aleph_0$
is handled using Theorem \ref{existenceofcg+}.

\begin{Theorem}\label{KjSh}[Kojman-Shelah \cite{KjSh 409}] Suppose that $\kappa$ and
$\lambda$ are regular such that $\lambda< 2^\kappa$.
Further suppose that either 
\begin{description}
\item{(a)} $\kappa=\aleph_0, \lambda=\aleph_1$ and there is
an $S_{\aleph_0}^{\aleph_1}$-club guessing
sequence; or 
\item{(b)} $\kappa^+ <\lambda$.
\end{description}
Then there is no universal
linear order of size $\lambda$, moreover the universality number of the
class of linear orders of size $\lambda$ is at least $2^\kappa$.
\end{Theorem}

\begin{Proof} (sketch) 
The proof uses the method of Construction and Preservation.
First let us fix 
a club guessing sequence $\bar{C}=\langle C_\delta:\,\delta\in S\subseteq
S^\lambda_\kappa\rangle$ and if $\kappa>\aleph_0$ assume also that
this sequence is chosen in conjunction with a sequence $\bar{P}$ to form a 
truly tight $(\kappa,\lambda)$-guessing sequence.
In particular we assume that the order type $i^\ast(\delta)$ of $C_\delta$ is always
$\kappa$. 
Recall the notation $\langle \alpha^\delta_i:\,i<\kappa\rangle$ for the
increasing enumeration of $C_\delta$. 

\begin{Lemma}\label{constr.lemma}[Construction Lemma] For every $A\subseteq \kappa$, there is
a linear order $L_A$ and its filtration
${\bar L}_A$ such that for a club $C$ of $\lambda$, we have 
\[
(\delta\in C\,\,\&\,\,C_\delta\subseteq\delta)\implies {\rm inv}_{\bar{L}_A,\bar{C}\textcolor{red}{,\delta}}(\delta)=A.
\]
\end{Lemma}

\begin{Proof} (sketch) The Construction Lemma has a simpler proof in the case (a) of
the theorem, which is the case we shall prove. Then we shall comment on the changes
needed to cover case (b). So assume that $\kappa=\aleph_0$, $\lambda=\aleph_1$ and
$\bar{C}$ is an $S^{\aleph_1}_{\aleph_0}$-club guessing sequence. Let $A\subseteq
\omega$ be given.

Recall that by a {\em cut}
of a linear order we mean an initial segment of the order, and we say that the cut is
{\em realised} if it has the least upper bound. If $L\subseteq L'$ are linear orders
and $D$ is a cut of $L$, then a cut $D'$ of $L'$ {\em extends} $D$ if $D'\cap L=D$.
When speaking of a linear order $(L,<_L)$ we may refer to cuts of $L$ or cuts of
$<_L$, as is more convenient for the context.
By $\eta$ we denote the order type of the rationals.

We shall defines the order $<_{L_A}$ on $\omega_1$ by inductively defining
a strictly increasing sequence $\langle \gamma_i:\,i<\omega_1\rangle$ of countable
ordinals and
defining $<_i\deq<_{L_A}\rest{\gamma_i}$ at the step $i$ of the induction. The requirements
of this induction will be:
\begin{description}
\item{(i)} for every $i<j$ there is a cut of $<_i$ realised in $\gamma_{j+1}$ and
not realised in $\gamma_{j}$;
\item{(ii)} if $i<j<k$ and $D$ is a cut of $\gamma_i$ realised in $\gamma_j$ but not in 
$\gamma_i$, then there is a cut of $D$ of $\gamma_k$ that extends $D$ and that is
realised in $\gamma_{k+1}$ but not in 
$\gamma_k$;
\item{(iii)} if $D$ is a cut of $\gamma_i$ realised in $\gamma_{i+1}$ and
not realised in $\gamma_{i}$, then the $<_{i+1}$-order type of $\{x\in \gamma_{i+1}:\,
x\mbox{ realises }D\}$ is $\eta$.
\end{description}
The starting point of the induction
is $\gamma_0=\omega$ where we let $<_0$ be $\omega$ ordered in the order type of the rationals.
At limit $i$ we define $\gamma_i$ to be the $\sup_{j<i}\gamma_j$
and the order is defined as the union of the orders constructed so far.

At the stage $\gamma_{i+1}$ we ask ourselves if
$i$ is a {\em good point} of the construction, which means that $i$ is
a limit ordinal and $\gamma_i=i$. If so we then ask if $C_\delta$
consists of good points. If the answer to both of these question is affirmative
we proceed to define a sequence $\langle D_n:\,n<\omega\rangle$
of cuts, such that each $D_n$ is a cut of $\alpha^i_{n}$ not
realised in $<_{\alpha^i_n}$ and $D_n\subseteq D_{n+1}$.
In addition we require that $D_n$ is realised in $<_{\alpha_{n+1}^i}$ iff $n\in A$.
If $A=\emptyset$ we let $D_0$ be a cut of $<_0$ not realised in $i$, which 
exists as there are $2^{\aleph_0}$ cuts of $<_0$ and $i$ is countable.
We let $D_n$ be the extension of $D$ to $\alpha^i_n$, so
$D_n\deq\{x<\alpha^i_n:\,(\exists y\in D)(x<_{\alpha^i_n} y)\}$. Otherwise,
let $n_0$ be the first element of $A$ and let $D_{n_0}$ be a cut of $\alpha^i_{n_0}$
that is realised in $\alpha^{i}_{n_0+1}$ and not in $\alpha^i_{n_0}$. For
$m<n_0$ let $D_m=D_{n_0}\cap \alpha^i_m$. If $A=\{n_0\}$, then since there is the
order type $\eta$ of elements of $\alpha^{i}_{n_0+1}$ that realise  $D_{n_0}$
(by requirement (ii)), there is a cut $D$ of $i$ that extends $D_{n_0}$ and
that is not realised in $i$. For $n> n_0$ we let $D_{n}$ be $D\cap \alpha^i_n$.
Otherwise, let $n_1=\min A\setminus (n_0+1)$ and let $D_{n_1}$ be a cut of
$\alpha^i_{n_1}$ extending $D_{n_0}$, which is realised in $\alpha_{n_1+1}^i$
and not in $\alpha_{n_1}$, which exists by (iii) above. Then we continue
similarly to the previous case. In any case, we have constructed the increasing
sequence of cuts as required, and letting $D^\ast$ be their union
we then let $i$ realise $D^\ast$ in $<_L\rest (i+1)$. We extend the order by
transitivity.

Now we still have to assure that the requirements of the induction are preserved,
which can be done by amalgamating countably many ordinals to $i+1$ in the way
requested by the requirements. The sup of all these ordinals is then defined
to be $\gamma_{i+1}$. If $i$ is not a good point we do not have to take special
care of $i$ but instead proceed just as in this paragraph.

At the end of the induction we let $<_{L_A}$ be the union of the orders constructed,
$\bar{L}_A=\langle <_i:\,i<\omega_1\rangle$ and $E$ a club of good points of the
construction. The construction was made so that at any $\delta\in E$ the required invariant
is guaranteed to be achieved only if $C_\delta\subseteq E$.

This finishes the proof of the Construction Lemma in the case $\kappa=\aleph_0$.
For larger $\kappa$ things become more complex, as one also has to handle the
limit points of cofinality $<\kappa$. This is a difficulty familiar from
classical constructions, such as that of a Suslin tree from a $\diamondsuit$,
where one in addition uses a $\square$ sequence at cardinals larger than $\aleph_1$.
In this case the construction can be carried through thanks to the square-like
properties of a truly tight guessing sequence. See \cite{KjSh 409} for details.
$\eop_{\ref{constr.lemma}}$
\end{Proof}

We also need the Preservation Lemma. 

\begin{Lemma}\label{preser.lemma}[Preservation Lemma]
Suppose that $L$ and $L'$ are linear orders
with universe $\lambda$ and with filtrations $\bar{L}$ and $\bar{L'}$ respectively, while
$f:\,L\into L'$ is an order-preserving injection. Then there
is a club $E$ of $\lambda$ such that for every $\delta\in S^\lambda_\kappa$
satisfying $C_\delta\subseteq E$, we have
\[
{\rm inv}_{\bar{L},\bar{C}\textcolor{red}{,\delta}}(\delta)={\rm inv}_{\bar{L'},\bar{C}\textcolor{red}{,\delta}}(f(\delta)).
\]
\end{Lemma}

\begin{Proof}(sketch) We start by defining a model $M$ with universe $\lambda$,
order relations $<_L$, $<_{L'}$ and $<$ (the ordinary order on the ordinals)
and the function $f$. Let $E$ be a club of
$\delta<\lambda$ such that $\delta\in E$ implies that $M\rest\delta\elementary M$
and the universe of both $L_\delta$ and $L'_\delta$ is $\delta$. \textcolor{red}{Note that then $f(\delta)\ge\delta$}.
Suppose that
$\delta\in E$ is such that $C_\delta\subseteq E$, and we shall prove that
${\rm inv}_{\bar{L},\bar{C}\textcolor{red}{,\delta}}(\delta)={\rm inv}_{\bar{L'},\bar{C}\textcolor{red}{,\delta}}(f(\delta))$. The
more difficult direction of the proof is the inclusion $\subseteq$. So suppose that
$i\in {\rm inv}_{\bar{L},\bar{C}\textcolor{red}{,\delta}}(\delta)$, hence there is $\beta\in L_{\alpha^\delta_{i+1}}
\setminus L_{\alpha^\delta_{i}}$ satisfying that 
\[
\{x<_L\beta:\,x\in L_{\alpha^\delta_{i}}\}=\{x<_L\delta:\,x\in L_{\alpha^\delta_{i}}\}.
\]
We would like to claim that $f(\textcolor{red}\beta)$ witnesses that $i\in {\rm inv}_{\bar{L'},\bar{C}\textcolor{red}{,\delta}}(f(\delta))$,
and it does follow from the choice of $E$ and the fact that $C_\delta\subseteq E$ that

\[
\{y<_{L'}f(\beta):\,y\in f``(L_{\alpha^\delta_{i}})\}=\{y<_{L'}f(\delta):\,y\in
f``(L_{\alpha^\delta_{i}})\}.
\]
However, the problem is that $f$ is not necessarily onto.
As $L'$ is a linear order we have $f(\beta)<_{L'} f(\delta)$ or $f(\delta)<_{L'} f(\beta)$
(\textcolor{red}{the} equality cannot occur by the choice of $E$).
Let us suppose that the former is true, the latter
case is symmetric. Suppose that $f(\beta)$ does not witness that
$i\in {\rm inv}_{\bar{L'},\bar{C}\textcolor{red}{,\delta}}(f(\delta))$, this then means that there is $\gamma \in
L'_{\alpha^\delta_i}$ (so $\gamma< \alpha^\delta_i$) 
such that $f(\beta)<_{L'}\gamma<_{L'} f(\delta)$. Observe that there is no $\varepsilon
\in L_{\alpha^\delta_i}$ such that $\gamma\le_{L'} f(\varepsilon)\le_{L'}f(\delta)$,
by the choice of $\beta$. 

Consider $T\deq\{x:
\,(\neg\exists q\in L_{\alpha^\delta_i})\,\gamma
<_{L'}f(q)<_{L'} x\}$. We claim that $T\cap\alpha^\delta_i$ is exactly the set
$\{\zeta<\alpha^\delta_i:\,\zeta<_{L'} f(\delta)\}$. Namely if
$\zeta<\alpha^\delta_i$ and $\zeta<_{L'} f(\delta)$ and $\zeta\notin T$ then there is
$\varepsilon\in L_{\alpha^\delta_i}$ such that $\gamma<_{L'}<f(\varepsilon)<_{L'}\zeta
<_{L'}f(\delta)$, a contradiction. On the other hand, if for some $\zeta<\alpha^\delta_i$
we have $f(\delta)<_{L'}\zeta$ then in $M$ it is true that there is $q$ such that
$\gamma<_{L'}f(q)<_{L'}\zeta$, as $\delta$ is such a $q$. By elementarity it is true
that there is such $q\in L_{\alpha^\delta_i}$, so $\zeta\notin T$. Hence we have shown
the existence of a $\xi$ such that $\{\zeta<\alpha_i^\delta:\,\zeta<_{L'}\xi\}$
is exactly 
\[
\{x<\alpha^\delta_i:
\,(\neg\exists q\in L_{\alpha^\delta_i})\,\gamma
<_{L'}f(q)<_{L'} x\}.
\]
By elementarity there must be such $\xi\in L'_{\alpha^\delta_{i+1}}\setminus
L'_{\alpha^\delta_{i}}$, and then this $\xi$ shows that $i\in 
{\rm inv}_{\bar{L'},\bar{C}\textcolor{red}{,\delta}}(f(\delta))$.
$\eop_{\ref{preser.lemma}}$
\end{Proof}

To finish the proof of the theorem, suppose that there were a family 
$\{L_i:\,i<i^\ast\}$ of linear orders of size $\lambda$
for some $i^\ast<2^\kappa$, such that every linear order of size $\lambda$
embeds into some $L_i$. We may assume that the universe of each
$L_i$ is $\lambda$. Let $\bar{L}_i$ be any filtration of $L_i$ and let ${\mathcal B}$
be the family of all $B\subseteq \kappa$ such that for some $x$, $i$ \textcolor{red}{and $\delta$} we have
${\rm inv}_{\bar{L_i},\bar{C} \textcolor{red}{,\delta}}(x)=B$. Then the size of $\BB$ is at most $\lambda\cdot
\card{i^\ast}$, which is $< 2^\kappa$. Hence there is by Construction Lemma a linear
order $L_A$ with universe $\lambda$ and its filtration $\bar{L}_A$
such that for a club $C$ of $\lambda$
we have $\delta\in C\implies {\rm inv}_{\bar{L}_A,\bar{C}\textcolor{red}{,\delta}}(\delta)=A$.
Suppose that $f:\,L^\ast\into L_i$ is an embedding and let $E$ be a club guaranteed to
exist by the Preservation Lemma. Let $\delta\in C$ be such that $C_\delta\subseteq E$.
Then ${\rm inv}_{\bar{L}_i,\bar{C}\textcolor{red}{,\delta} }(\delta)(f(\delta))=A$, a contradiction with the choice of $A$.
$\eop_{\ref{KjSh}}$
\end{Proof}

\section{A few more words
on orders and orderable structures}\label{extras}
The Kojman-Shelah method
can be ramified to give results on cardinals $\lambda$ that are not necessarily
regular. Using another pcf staple, 
the covering number ${\rm cov}(\lambda,\mu,\theta,\sigma)$, they were able to strap 
together the negative universality results on the regular cardinals below a given
singular cardinal to obtain a negative universality result for the singular. Note that
by classical results about special models (see \cite{ChKe}) there is a universal
linear order (or any other first order theory of a sufficiently
small size) in any strong limit uncountable cardinality. The question then becomes what
happens if for example $\aleph_\omega$ is not a strong limit. The answer is that then
there is no universal linear order then. Specifically:

\begin{Theorem}\label{KjShsing}[Kojman-Shelah \cite{KjSh 409}] Suppose that $\lambda$ is 
a singular cardinal which is not a strong limit and it satisfies that either
\begin{description}
\item{(a)} $\aleph_\lambda>\lambda$ or
\item{(b)} $\aleph_\lambda=\lambda$ but $\card{\{\mu<\lambda:\,\aleph_\mu=\mu\}}<
\lambda$ and either $\cf(\lambda)=\aleph_0$ or $2^{<\cf(\mu)}<\mu$,
\end{description}
then there is no universal linear order of size $\lambda$.
\end{Theorem}

Linear orders are representatives of theories $T$ that have the {\em strict order property},
which means that there is a formula $\varphi(\bar{x};\bar{y})$ such that in the
monster model $\CC$ of $T$ there are $\bar{a}_n$ for $n<\omega$ such that
for any $m,n<\omega$
\[
\CC\models ``(\exists \bar{x})[\neg\varphi(\bar{x};\bar{a}_m)\wedge \varphi(\bar{x};\bar{a}_k)]"
\mbox{ iff }m<n.
\]
Other examples of first order theories that have the strict order property are
Boolean algebras, partial orders, lattices, ordered
fields, ordered groups and any unstable complete theory that does not have the independence
property (see \cite{Sh -c}). Using the fact that the strict order property of $T$ allows for
coding of orders into models of $T$ 
and that there is a quantifier-free definable order in
the above (non-complete) theories, \S 5 of \cite{KjSh 409} shows that the existence of 
a universal element in any of these theories at a cardinal $\lambda$ implies the existence
of a universal linear order of size $\lambda$. Therefore the negative universality
results stated above also apply to these theories. 

A different approach to drawing
conclusions about the universality problem in one class knowing the behaviour of another
class is taken by Katherine Thompson
\cite{KTh} and \cite{KTh2}, who uses functors that preserve the embedding structure. She reproves
the Kojman-Shelah conclusion about universality of partial orders versus
that of linear orders and connects certain classes of graphs with certain classes
of strict orders. This approach is useful also when one moves
from the first order context, for example to the class of orders
that omit chains of a certain type. The simplest case are orders that omit
infinite descending sequence. Universality
is resolved trivially in the class of well orders, as follows from the example of
the ordinals, but by changing the context
to that of well-founded partial orders with some extra requirements one obtains a
different situation. This type of problem is the subject of \cite{KTh}.

\section{Model theory}\label{models} There is
a model-theoretic motivation behind a 
an attempt to deliver a general method of approach to the universality problem,
stemming from Shelah's programme of classification theory. This is very well 
described in \S5 of Shelah's paper on open questions in model theory,
\cite{Sh 702}.
Namely it may be hoped that the behaviour of a theory with respect to the universality
would classify the theory as `good' if it 
can admit a small number of universal models even when the relevant instances of $GCH$ 
fail, while a bad theory would rule out small universal families as soon as
$GCH$ would be sufficiently violated (recall that the situation in the presence of
$GCH$ is information-free here, as all first order countable theories e.g. have a universal
element in every uncountable cardinal (see \cite{ChKe}). Such a division would be used to 
classify unstable theories, with the hoped for result similar to the classification
of stable versus unstable theories where a model-theoretic property of stability of
a countable theory
was closely connected with the number of of nonisomorphic models a theory may 
have at an uncountable cardinal, through the celebrated Shelah's Main Gap Theorem
(see \cite{Sh -c}).
The idea of using universality in a similar manner has proved to be quite
successful, and although no precise model-theoretic equivalent has been found as of yet,
there is much information available about the existing model-theoretic properties. One can
find a rather detailed description of the present state of knowledge in \cite{DjSh 710}
where there is also a precise
definition of the proposed division from the set-theoretic point of view, that is
what is meant by being 'good' (referred to as {\em amenable}) and `very-bad'
({\em highly non-amenable}) from the universality point of view. In this
paper we mostly concentrate on the highly non-amenable theories, which can
be defined by

\begin{Definition}\label{hna}
A theory $T$ is said to be {\em highly non-amenable} iff
for every large enough regular cardinal $\lambda$ and $\kappa<\lambda$ such that
there is a truly tight $(\kappa,\lambda)$ club guessing sequence
$\langle C_\delta:\,\delta\in S\rangle$
the smallest number of models of $T$ of size $\lambda$ needed 
to embed all models of $T$ of that size is at least $2^\kappa$. $T$ is
{\em highly non-amenable up to $\kappa^\ast$} if the above characterisation
is not necessarily true, but it is true whenever $\kappa<\kappa^\ast$.
\end{Definition}

In model theory one usually works with complete theories, while our examples
above were not necessarily so (for example we worked with the theory of a linear
order). We adopt the convention that when speaking of a complete theory by
an embedding we mean an elementary embedding, and otherwise we just mean an
ordinary embedding. With this clause Definition \ref{hna} makes sense in both
contexts, and
the work of Kojman and Shelah presented in \S\ref{orders} showed that linear orders
and theories with the strict order property are highly non-amenable. Shortly after
this work the same authors in \cite{KjSh 447} proceeded to show (Theorem 4.1 + Theorem 5.1
of \cite{KjSh 447}):

\begin{Theorem}\label{447} [Kojman-Shelah \cite{KjSh 447}]
Countable stable unsuperstable theories are
highly non-amenable up to $\aleph_1$. In general stable unsuperstable theories $T$
are highly non-amenable up to their stability cardinal $\kappa(T)$.
\end{Theorem}

Here we use the notion of the stability cardinal $\kappa(T)$ defined as the
minimal cardinal $\kappa$ such that for every set $A\subseteq{\mathfrak C}$ and
a type $p$ over $A$ there is $B\subseteq A$ such that $\card{B}<\kappa$ and
$p$ does not fork over $B$. It is proved in Shelah's Stability Spectrum Theorem
\cite{Sh -c} that for any stable $T$ we have $\kappa(T)\le\card{T}^+$ and for every
$\lambda$ we have $T$ is stable in $\lambda$ iff $\lambda=\lambda^{<\kappa(T)}$
and either $\lambda\ge 2^{\aleph_0}$ or $\lambda\ge$ the number $D(T)$ of parameter-free
types of $T$ in ${\mathfrak C}$. A countable complete first order theory $T$ is stable
unsuperstable iff $\kappa(T)=\aleph_1$. In Theorem \ref{447},
as well as in many other applications of the
method, a major issue is how to define an invariant. Suppose that $T$ is, for
simplicity, a complete countable stable unsuperstable theory.

\begin{Definition}\label{withforks} Let $\lambda$ be regular, $\bar{C}=
\langle C_\delta:\,\delta\in S\rangle$ a 
$(\kappa,\lambda)$ tight club guessing sequence,  and
$N$ a model of $T$ of size $\lambda$ given with a continuous increasing filtration
$\bar{N}=\langle N_i:\,i<\lambda\rangle$. We define for a 
$\delta\in S$ and a tuple $\bar{a}$ of $N$,
\[
{\rm inv}_{\bar{N},\bar{C}}(\bar{a})\deq\{i<\kappa:\,\mbox{ the type of }
\bar{a}\mbox{ over }N_{\alpha^\delta_i}\mbox{ forks over }N_{\alpha^\delta_{i+1}}\}.
\]
\end{Definition}

As before, we have used the notation
$\langle\alpha^\delta_i:\,i<\kappa\rangle$ for the increasing enumeration of $C_\delta$. 
In the case of $\kappa=\aleph_0$ it suffices to deal with ordinary club guessing
sequences. We do not have the space to introduce the notion of forking here, but
the intuition behind Definition \ref{withforks} is similar to the idea behind the
invariant for linear orders: $i$ is in the invariant iff `something new happens'
at the stage $\alpha^\delta_{i+1}$, something that `reflects' the behaviour
of $\bar{a}$ with respect to $\bigcup_{i<\kappa} N_{\alpha^\delta_i}$. Recalling that
$\kappa<\kappa(T)$ is assumed may give a hint of how the (rather complex) proof
of \cite{KjSh 447} proceeds.

One can also try to what is meant by the `good' universality behaviour, and in
\cite{DjSh 710} we have tried to capture this using the notion of
{\em amenability}.

\begin{Definition}\label{amenability}
A theory $T$ is {\em amenable} iff whenever $\lambda$ is an uncountable
cardinal larger than the size of $T$ and satisfying $\lambda^{<\lambda}=\lambda$
and $2^\lambda=\lambda^+$, while $\theta$ satisfies $\cf(\theta)>\lambda^+$,
there is a $\lambda^+$-cc $(<\lambda)$-closed forcing notion that forces $2^\lambda$
to be $\theta$ and assures that in the extension there is a family
$\FF$ of $<\theta$ models of $T$
of size $\lambda^+$ such that every model of $T$ of size $\lambda^+$ embeds into
one of the models in $\FF$. Localising at a specific $\lambda$ we obtain the
definition of amenability at $\lambda$.
\end{Definition}

The point is that no theory can be both amenable and highly non-amenable. Namely
suppose that a theory $T$ is both amenable and highly non-amenable, and let $\lambda$ be
a large enough regular cardinal while $V=L$ or simply
$\lambda^{<\lambda}=\lambda$ and $\diamondsuit(S^{\lambda^+}_\lambda)$
holds.
Let $P$ be the forcing exemplifying that $T$
is amenable. Clearly there is a truly tight $(\lambda,\lambda^+)$
club guessing sequence $\bar{C}$ in $V$, and since the forcing $P$ is $\lambda^+$-cc,
every club of $\lambda^+$ in $V^P$ contains a club of $\lambda^+$ in $V$,
hence $\bar{C}$ continues to be a truly tight $(\lambda,\lambda^+)$
club guessing sequence in $V^P$. Then on the one hand we have that in $V^P$, the universality
number of models of $T$ of size $\lambda$,
univ$(T,\lambda^+)$, is at least $2^\lambda$ by the high non-amenability, while univ$(T,\lambda^+)<
2^\lambda$ by the choice of $P$, a contradiction.
In \cite{DjSh 614}, building on the earlier work of Shelah in
\cite{Sh 457}, we gave an axiomatisation of elementary classes 
that guarantees that the underlying theory is amenable.
Shelah proved in \cite{Sh 500} that all countable simple theories are
amenable at all successors of regular $\kappa$ satisfying
$\kappa^{<\kappa}=\kappa$. (Note that even though all simple theories are stable, this is
not in contradiction with Theorem \ref{447},
as there it is only proved that countable stable unsuperstable theories
are highly non-amenable up to $\aleph_1$).
In that same paper Shelah introduced a hierarchy
of complexity for first order theories, and showed that high non-amenability appears
as soon as a certain level on that hierarchy is passed.
Details of this
hierarchy are given in the following definition:

\begin{Definition}\label{nsops} Let $n\ge 3$ be a natural
number.
A formula $\varphi(\bar{x},\bar{y})$
is said to exemplify the $n$-{\em strong order
property} of $T$,
$SOP_n$, if $\llg(\bar{x})=\llg(\bar{y})$, and there are
$\bar{a}_k$  for $k<\omega$, each of length ${\llg(\bar{x})}$
such that
\begin{description}
\item{(a)} $\models\varphi[\bar{a}_k,\bar{a}_m]$ for $k<m<\omega$,

\item{(b)} $\models \neg(\exists\bar{x}_{0},\ldots,\bar{x}_{n-1})[\bigwedge\{
\varphi(\bar{x}_{\ell},\bar{x}_{k}):\,{\ell},k<n\mbox{ and }k={\ell}+1\mod n\}].$
\end{description}
\end{Definition} 

The following were proved in \cite{Sh 500}: the hierarchy above describes
a sequence $SOP_n\,(3\le n<\omega)$ of properties of strictly increasing
strength such that the theory of a dense linear order possesses all the
properties, while on the other hand no simple theory can have
the weakest among them, $SOP_3$. The property
$SOP_4$ of a theory $T$ implies that $T$ 
is highly non-amenable. In the light of these results it
might then be asked
if $SOP_4$ is a characterisation of high non-amenability,
A partial solution appears in
\cite{DjSh 710}. There we considered a property of theories that we called {\em oak property},
as its prototypical example is a tree of the form ${}^{\kappa\ge}\lambda$ equipped
with restriction where we can express that $\eta\rest\alpha=\nu$ for $\eta\in 
{}^{\kappa}\lambda$, $\alpha<\kappa$ and $\nu \in {}^{\kappa >}\lambda$. This property
is also a generalisation of the theory of infinitely many independent equivalence
relations $T^\ast_{\rm feq}$, see \cite{DjSh 710}.
The formal definition
is:

\begin{Definition}\label{Prs}
A theory $T$ is said to {\em satisfy}
the {\em oak property as exhibited by
a formula} $\varphi(\bar{x},\bar{y}, \bar{z})$ iff
for any infinite $\lambda,\kappa$ there are $\bar{b}_\eta
(\eta\in{}^{\kappa>}\lambda)$, $\bar{c}_\nu (\nu \in
{}^{\kappa}\lambda)$ and $\bar{a}_i (i<\kappa)$
such that
\begin{description}
\item{(a)} $[\eta\initial\nu\,\,\&\,\,\nu \in {}^\kappa\lambda]
\implies\varphi[ \bar{a}_{\llg(\eta)},\bar{b}_\eta, \bar{c}_\nu]$,
\item{(b)} If 
$\eta\in {}^{\kappa>}\lambda$ and $\eta\,\hat{}\,\langle\alpha\rangle\initial
\nu_1\in {}^\kappa\lambda$ and $\eta\,\hat{}\,\langle\beta\rangle\initial
\nu_2\in {}^\kappa\lambda$, while $\alpha\neq\beta$ and $i>\llg(\eta)$,
\underline{then} $\neg\exists \bar{y}\,
[\varphi(\bar{a}_i, \bar{y}, \bar{c}_{\nu_1})
\wedge \varphi(\bar{a}_i, \bar{y},\bar{c}_{\nu_2} )]$,
\end{description}
and in addition $\varphi$ satisfies
\begin{description}
\item{(c) $\varphi(\bar{x},\bar{y}_1, \bar{z})\wedge
\varphi(\bar{x},\bar{y}_2, \bar{z})\implies \bar{y}_1=\bar{y}_2$.
}
\end{description}
\end{Definition}

Shelah proved in \cite{Sh 457} that $T^\ast_{\rm feq}$
exhibits a non-amenability behaviour
provided that some cardinal arithmetic assumptions close to the failure of the
singular cardinal hypothesis are satisfied. This does not
necessarily  imply high non-amenability as it was proved also in \cite{Sh 457}
that this theory is in fact amenable at any cardinal which is the successor of
a cardinal $\kappa$ satisfying $\kappa^{<\kappa}=\kappa$.
In \cite{DjSh 710} we generalised the first of these two results by
showing that any theory with oak property satisfies the same non-amenability
results as those of $T^\ast_{\rm feq}$, and we gave some more circumstances,
given in terms of pcf theory for when such non-amenability results hold.
The oak property cannot be made
a part of the $SOP_n$ hierarchy,
as \cite{DjSh 710} gave a theory which has oak, and is $NSOP_3$, while the model
completion of the theory of triangle free graphs is an example of a $SOP_3$
theory which does not satisfy the oak property.
On the other hand it is also proved in \cite{DjSh 710} that no oak theory
is simple. Further considerations of the oak property appear in \cite{Sh 820}
where it is proved that (under an interpretation of
what it means for a class to have oak) that the class of groups has this
property. That paper also gives further universality results in the context of
abelian groups.

\section{Some applications in analysis and topology}\label{analysis} There is
a rich literature concerning the universality problem
in the various classes of compact spaces coming from analysis, such as Corson
and Eberlein compacta, where by an embedding we usually mean
the existence of a continuos surjection, see e.g. \cite{AB}.
Many of these questions very nicely resolved by the $\sigma$-functor of
Todor\v cevi\'c, \cite{Todor}, which gives for every such space $K$
another space $\sigma(K)$ in the same class such that $\sigma(K)$ is not a
continuos image of $K$. The class of uniform Eberlein compacta, which are those
compact spaces that are homeomorphic to a
weakly compact subspace of a Hilbert space, seem to be an odd one out in these
problems, since neither the method of generalised Szlenk invariants as employed
in \cite{AB} not the $\sigma$ functor give any results in this class. The reason
is, as the authors of \cite{AB} observed, that their invariant defined
for all Eberlein spaces, trivialises in the case that the Eberlein is uniform,
while $\sigma(K)$ for a uniform Eberlein compact space $K$
is an Eberlein compact but not necessarily uniform. Murray Bell 
in \cite{Bell} made a major advance in
the universality problem of the uniform Eberlein compacta, which had been
completely open since the 1977 paper  \cite{BeRuWa} of Yoav Benyamini, Mary Ellen Rudin and
Michael Wage posed it. Namely, Bell defined a certain algebraic structure, the so
called c-algebra, and he proved that there is a universal UEC of weight $\lambda$ iff
there is a universal c-algebra of size $\lambda$. In the same paper
Bell showed that if $2^{<\lambda}=
\lambda$, there is a c-algebra
of size $\lambda$ which is universal not just under ordinary embeddings,
but also under a stronger notion of a c-embedding. We shall call such 
algebras c-{\em universal}, definitions follow. He also provided
negative consistency results in models obtained by adding Cohen subsets to a regular
cardinal.

\begin{Definition}\label{embedding}(1) 
A subset $C$ of a
Boolean algebra $B$ has {\em the nice property}
if for no finite $F\subseteq C$ do we have $\bigvee F=1$.
A Boolean algebra $B$ is a c-{\em algebra} iff there is
a family $\langle A_n:\,n<\omega\rangle$ of pairwise disjoint subsets of $B$
each consisting of pairwise disjoint elements, whose union
has the nice property and generates $B$. 

{\noindent (2)} If $B_l$ for $l\in \{0,1\}$ are c-algebras with fixed
sequences $\langle A_n^l:\,n<\omega\rangle$ of subsets exemplifying that $B_l$ is
a c-algebra,
then a 1-1 Boolean 
homomorphism $f:\,B_0\into B^\ast$ is a c-{\em embedding} iff
$f``A_n^l\subseteq A_n^l$ for all $n<\omega$. 
\end{Definition}

Note that the notion of a c-algebra is not first order so the Kojman-Shelah results
from \cite{KjSh 409} do not directly apply.
We showed in \cite{Dz}
that for no regular cardinal $\lambda>\aleph_1$
with $2^{\aleph_0}>\lambda$ can there exist $<2^{\aleph_0}$
c-algebras of size $\lambda$ such that every
c-algebra of size $\lambda$ embeds into one of them. These results were continued in \cite{Dz2}
which contains both negative and positive results about the existence
of universal c-algebras and UEC. 
On the other hand, we proved a
positive consistency result showing that under certain non-GCH
assumptions there can be a family of UEC of a relatively
small size ($\lambda^{++}< 2^{\lambda^+}$) each of which has weight
$\lambda^+$ and which are jointly universal for UEC of weight
$\lambda^+$. The negative results are the ones relevant to this paper, and they
were obtained using the method of Kojman-Shelah invariants. The appropriate definition
in this context turned out to be the following:

\begin{Definition} Let $\lambda$ be a regular cardinal and 
$\langle C_\delta:\,\delta\in S\rangle$ a club guessing sequence on $\lambda$,
with $C_\delta=\langle \alpha^\delta_i:\,i<i^\ast\rangle$ an increasing enumeration.
Let $B$ be a c-algebra of size $\lambda$ with
a filtration $\bar{B}$, and we assume that $\langle A_n:\, n<\omega\rangle$
is a fixed sequence demonstrating that $B$ is a c-algebra.
Suppose that $\delta\in S$ and define
for $\delta\in S$ and
$b\in B\setminus B_\delta$ 
\[
{\rm inv}_{\bar{B}, \bar{C}}(b)\deq\left\{
i<i^\ast:(\exists m\ge 1)
(\exists y\in A_m\cap B_{\alpha^\delta{i+1}}\setminus B_{\alpha^\delta_i})
\,[y\ge b]\right\}.
\]
\end{Definition}

It is interesting to note that the mutual exclusiveness of amenability and
high non-amenability as defined in \S\ref{models} does not apply if these
definitions are taken in their obvious translation to the non-first order
context, and that the class of c-algebras exemplifies that.

A recent application of the method of invariants comes from Kojman-Shelah's work on
almost isometric embeddings between metric spaces, \cite{KjSh 827}. A map $f:\,X\into Y$
between metric spaces is said to be {\em Lipshitz with constant} $r>0$ if for
every $x,y\in X$ we have $d_Y(f(x),f(y))<r\cdot d_X(x,y)$. $X$ is {\em almost isometrically
embeddable into} $Y$ iff for every $r>1$ there is a continuous injection
$f:\,X\into Y$ such that both $f$ and $f^{-1}$ are Lipshitz with constant $r$,
which is called {\em bi-Lipshitz with constant} $r$.
Among many interesting results about such embeddings, \cite{KjSh 827} also
gives 

\begin{Theorem}\label{metricth}[Kojman-Shelah, \cite{KjSh 827}]
If $\aleph_1< \lambda<2^{\aleph_0}$ 
is regular then for every $\kappa<2 ^{\aleph_0}$ and metric spaces
$\{(X , d_i) : i < \kappa\}$ of size $\lambda$,
there exists a metric space of size $\lambda$ that is not almost isometrically
embeddable into any $( X, d_i)$.
\end{Theorem}

The proof of the theorem again uses the method of invariants, but with a twist.
Namely one defines two kind of invariants, ${\rm inv}^{\rm dom}$ and ${\rm inv}^{\rm rng}$,
as follows, where we are using the same notation for club guessing sequences as above:

\begin{Definition} Suppose that $(X,d)$ is a metric space with universe $\lambda$,
$\bar{C}$ is an $S_{\aleph_0}^\lambda$ club guessing sequence, $\delta\in S^\lambda_{\aleph_0}$,
$\beta>\delta$ and $K\ge 1$ is an integer.  We consider $X$ as being given in the
filtration $\bar{X}=\{\alpha:\,\alpha<\lambda\}$. Then
\[
{\rm inv}^{\rm dom}_{\bar{C},X,\delta}(\beta)=
\{n<\omega:\,d(\beta,\alpha^\delta_n)/d(\beta,\alpha^\delta_{n+1})>2 K^2\}
\mbox{ and }
\]
\[
{\rm inv}^{\rm ran}_{\bar{C},X,\delta}(\beta)=
\{n<\omega:\,d(\beta,\alpha^\delta_n)/d(\beta,\alpha^\delta_{n+1})>4 K^4\}.
\]
\end{Definition}

The Preservation Lemma then says in particular
that if $f:\,X\into Y$, where both $X$ and
$Y$ are metric spaces with universe $\lambda$, is bi-Lipshitz with constant $K$, then
there is a club $E$ of $\lambda$ such that for every $\delta\in E\cap S^\lambda_{\aleph_0}$
and $\beta>\delta$, we have $f(\beta)>\delta$ and ${\rm inv}^{\rm dom}_{\bar{C},X,\delta}(\beta)=
{\rm inv}^{\rm rng}_{\bar{C},Y,\delta}(f(\beta))$. Theorem \ref{metricth} is to be contrasted
with another theorem from \cite{KjSh 827}, which says that for any regular cardinal
$\lambda$ it is consistent that $2^{\aleph_0}>\lambda^+$ and there
$\lambda^ +$ separable metric spaces of size $\lambda$ such that every separable
metric space of size $\lambda$ almost-isometrically embeds into one of them. 
Earlier results about universality of metric spaces under different kinds of embeddings
and involving the method of invariants were obtained by Shelah in \cite{Sh 552}.

Model theory of metric spaces is also one of the subjects of 
Alex Ustvyasov's Ph.D thesis \cite{Us} and his joint work with Shelah in
\cite{ShUs 837}, where they concentrate on complete metric spaces. A model-theoretic
approach to Banach spaces pays off in Shelah-Ustvyasov's paper \cite{ShUs 789}
where they prove that the appropriately axiomatised theory of Banach spaces
has $SOP_n$ for all $n\ge 3$, and hence draw the negative universality
results provided by $SOP_4$ (see \S\ref{models}), where the notion of embedding is isometry.
Note that if $\lambda=2^{<\lambda}>\aleph_0$ then there is an isometrically universal
Banach space of size $\lambda$. Universality results in Banach
spaces are quite well studied classically, maybe the most well known result in this vein
is that of Szlenk in \cite{Sz} who proved that there is no universal reflexive separable
Banach space.

\section{Some applications in algebra}\label{algebra} A very fruitful application of
the Kojman-Shelah method of invariants has been in the theory of infinite abelian groups,
which we shall take in their additive notation.
In \cite{KjSh 455} Kojman and Shelah study the problem of universality in several kinds of
groups under various kinds of embeddings. 
Many classes of groups simply have a universal element under ordinary embeddings
in every infinite cardinality, namely there is always
a universal group, universal $p$-group (for any prime $p$),
universal torsion group and universal torsion-free group (see \cite{KjSh 455}).
On the other hand there is
no universal reduced $p$-group. The situation becomes different when one restricts
the kind of embeddings and the kind of groups one considers. Of particular interest are 
{\em pure} embeddings, where a group monomorphism $f:\,H\into G$ is a pure embedding if
$f``H$ satisfies that for all $n<\omega$, $nf``H=NG\cap f``H$. In other words, $f``H$
is a {\em pure} subgroup of $G$.

The appropriate notion of the invariant here is

\begin{Definition}\label{invpure} 
Suppose that $\lambda>\aleph_0$ is regular cardinal, $G$ is an abelian group of
size $\lambda$ given with its filtration $\bar{G}$ and $\langle C_\delta:\,\delta\in S\subseteq
\lambda\rangle$ is a club guessing sequence on $\lambda$ where for each $\delta$ the increasing
enumeration of $C_\delta$ is $\langle \alpha^\delta_i:\,i<i^\ast_\delta\rangle$. For $g\in G$
and $\delta\in S$ we define
\[
{\rm inv}_{\bar{G},\bar{C},\delta}(g)\deq\{i<i^\ast_\delta:\,g\in \bigcup_{n<\omega}
((G_{\alpha^\delta_{i+1}}+nG)\setminus (G_{\alpha^\delta_i}+nG))\}.
\]
\end{Definition}

A Preservation Lemma can be proved for this type of invariants and pure embeddings. The paper
gives a number of constructions of various types of groups with the prescribed invariant,
which allows for the proof of several theorems, a selection of which is:

\begin{Theorem}\label{groups1} [Kojman-Shelah, \cite{KjSh 455}] Suppose that
$\lambda$ is regular and for some $\mu$, $\mu^+<\lambda<2^{\mu}$, while $p$ is any prime. Then
there is no 
\begin{description}
\item{(a)} purely universal separable $p$-group of size $\lambda$;
\item{(b)} universal reduced slender group of size $\lambda$;
\item{(c)} universal reduced torsion-free group of size $\lambda$.
\end{description}
\end{Theorem}

Research on the universality in various classes of groups was continued by Shelah in
\cite{Sh 456}, \cite{Sh 552} and \cite{Sh 622}, where he considered various classes
of groups under ordinary embeddings (so they are not assumed to be pure).
In \cite{Sh 456} the class considered is that of
$(<\lambda)$-{\em stable} abelian groups, which means that for every subset $A$ of $G$
of size $<\lambda$ the closure in $G$ of the subgroup $\langle A\rangle_G$
generated by $A$, defined as 
\[
{\rm cl}_G(\langle A\rangle_G)=
\{x:\,\inf_{y\in \langle A\rangle_G}(\min_{i>1}
\{2^{-i}:\,x-y\mbox{ divisible by }\Pi_{1<j<i}n_j\})=0\}
\]
for some conveniently chosen and fixed increasing sequence $\langle n_i:\,i<\omega\rangle$
of natural numbers $>1$.
This notion in particular includes strongly $\lambda$-free groups and it can be handled
using the same definition of invariant as in Definition \ref{invpure}. In \cite{Sh 552}
there is a deep analysis of how necessary this is, and it proceeds through a series
of results about classes of trees with $\omega+1$ levels, with the thesis that these are a
prototype for various classes of groups (deriving also some
surprising results as to what kind of trees one needs to look at here).
Of particular interest in \cite{Sh 552} are reduced
torsion free groups and reduced separable
abelian $p$-groups, but the paper indeed gives a very rich selection of
results on various classes
of both groups and trees. This research was continued in \cite{Sh 622} and 
\cite{Sh F319}, and as a combined result of, one has almost a complete calculation
of the universality spectrum of the reduced torsion free abelian groups and reduced 
separable $p$-groups, for example:

\begin{Theorem}\label{combination} [Shelah] Let $\lambda$ be an
infinite cardinal and $\KK_\lambda$ the class of reduced torsion
free abelian groups of size $\lambda$ considered under ordinary embeddings.
\begin{description}
\item{(a)} If $\lambda=\lambda^{\aleph_0}$ or $\lambda$ is singular of countable cofinality and
$(\forall\theta<\lambda)\theta^{\aleph_0}<\theta$, then there is a universal member of 
$\KK_\lambda$.
\item{(b)} If $\lambda<2^{\aleph_0}$, or for some
$\mu$ we have $\beth_\omega+\mu^+<\lambda=\cf(\lambda)<\mu^{\aleph_0}$,
then there is no universal member of 
$\KK_\lambda$.
\end{description}
\end{Theorem}

Some of the remaining cases of possible cardinal arithmetic assumptions were reduced to
some weak pcf assumptions the consistency of whose failure is not known. In \cite{Sh 552}
there are also results about modules. Model theoretic properties of groups 
related to universality are
studied in \cite{ShUs 789}, where it is proved that if $\GG$ is the ``universal
domain" (a monster model for groups) then it has $SOP_3$ and, surprisingly, it does not
have $SOP_4$.

\section{Some applications in graph theory and a representation theorem}\label{graphs}
Universality problem in the class of graphs has a particularly long tradition,
see for example the well known
Richard Rado's paper \cite{Ra} with the construction of the Rado graph. If one 
considers graphs with the ordinary notion of embedding (so the edges are kept, but 
not necessarily non-edges, also called weak embedding), then under $GCH$ there is a universal
object in every infinite cardinality, as follows for uncountable cardinalities
from the classical first order model theory and was also proved independently by Rado.
Similar results hold for the class of graphs omitting the complete
graph $K_n$ of size $n$, where $n\ge 3$.
It is a very interesting result of Shelah that even when $CH$ fails there
can be a universal graph of size $\aleph_1$, see \cite{Sh 175}, \cite{Sh 175a}.
Results in \cite{Sh 500} and \cite{DjSh 614} imply that the theory of graphs
is amenable. The situation becomes different when one restricts to graphs that
omit a certain structure. For example, the model completion of the
theory of triangle-free graphs is amenable
but the model completion of the theory of directed
graphs omitting directed cycles of length $\le 4$ has $SOP_4$ and is
hence highly non-amenable (see \cite{Sh 500}, where a number of other similar results
is given). Passing to graphs that omit an infinite structure, so exiting the realm
of the first order theories, the situation immediately becomes very different.
For example, it is a mathematical folklore (see \cite{KoSh 492}) that there is no
universal $K_{\aleph_0}$-free graph in any cardinality. Komj\'ath and Shelah in \cite{KoSh 492}
investigate the class of $K_{\kappa}$-free graphs and show that under $GCH$ a universal
exists in $\lambda$ iff $\kappa$ is finite or $\cf(\kappa)>\cf(\lambda)$. They also
give consistency results showing how much the universality number of this class can be when
it is known that there is no one universal element. There is a very rich literature
available on the problem of the existence of a universal member in various
classes of graphs, for example there is a complete classification
of countable homogeneous directed graphs and countable homogeneous
$n$-tournaments, obtained by Gregory Cherlin in the memoir \cite{Cherlin}.
We cannot even begin to do justice to this rich literature
in this survey, so we shall simply concentrate on the impact the club guessing
method has had. This will also give us an opportune way of closing this paper by a theorem
which very elegantly shifts the method of invariants from an arbitrary class of models
of size $<2^{\aleph_0}$ to a consideration of the structure of the subsets of
the reals, namely a Representation Theorem by Kojman.

A {\em ray} in a graph is a 1-way infinite path. A {\em tail} of such a ray is any infinite
connected subgraph, and two rays are {\em tail equivalent} if they have a common ray.
Consider the class $\KK$ of graphs $G$ that satisfy that for every vertex $v$ of $G$
the induced subgraph of $G$ spanned by $v$ has at most one ray, up to tail equivalence.
This class can also be described in terms of forbidding certain structures. Among
other theorems that Kojman proves about this class in \cite{Kj} is that for a
regular uncountable $\lambda$ the smallest size of a family of graphs in $\KK$ of
size $\lambda$ (denoted by $\KK_\lambda$)
needed to embed such graphs is at least $2^{\aleph_0}$. This theorem
actually follows from a representation obtained in the following

\begin{Theorem}\label{representation} 
[Kojman's Representation Theorem, \cite{Kj}] If $\lambda>\aleph_1$ is regular then there is
a surjective homomorphism from the structure $\KK_\lambda$ partially ordered by 
the embedding relation, to the structure $[{\Bbf R}]^{\le\lambda}$ partially ordered
by the subset relation.
\end{Theorem}

The method of invariants is still used here, where the Construction Lemma corresponds
to proving that the proposed map is surjective and the Preservation Lemma corresponds
to showing that the map is a homomorphism. The Representation Theorem has some advantages
over the Construction and Preservation approach because it allows for a smooth way to handle
singular cardinals. In an upcoming paper \cite{DzTh} we have used this method to
consider well-founded partial orders under rank-preserving embeddings and some other
classes, and to prove negative universality results analogous to those in \cite{Kj}.

Let us finish by mentioning
that guessing sequences stronger than club guessing are used in a recent paper of Shelah
\cite{Sh 706} to obtain negative results about
the universality of the class of graphs that omit complete bipartite graphs. The paper
also gives a complete characterisation of the universality problem in this class
under $GCH$.

\section{Some questions}\label{questions}
There are many open questions in this subject, and considering instances of universality
in a specific class is an interesting pursuit per se. We have selected two more
general questions that are in our view very important. The first is in model theory:

\begin{Question} Does $SOP_4$ characterise high non-amenability,
in other words does every highly non-amenable theory
has the $SOP_4$ property?
\end{Question}

The interest of this question is described in \S\ref{models}.
The second question is in set theory and calls for a finer understanding of our
forcing iteration techniques. Namely, the reader may have noticed that the definition
\ref{amenability} does not refer to the existence of universal family $\FF$ of size 1, 
namely the
universal model. The reason is that all we know how to do, in the generality of
the axioms of \cite{DjSh 614} or in specific forcing proofs of universality such
as \cite{MeVa} (where Alan Mekler and Jouko V\"a\"an\"anen produced consistently with $CH$
a family of $\aleph_2$ trees of size $\aleph_1$ with no uncountable branches and
universal under reductions) or the Kojman-Shelah's theorem about
separable metric spaces mentioned above, is to produce $\lambda^{+}$ models
of size $\lambda$ jointly universal for models of size $\lambda$. Hence the
question is:

\begin{Question} Suppose that
a theory $T$ is amenable, $\lambda$ is an uncountable
cardinal larger than the size of $T$ and satisfying $\lambda^{<\lambda}=\lambda$
and $2^\lambda=\lambda^+$, while $\theta$ satisfies $\cf(\theta)>\lambda^+$. Can
one find a cardinality preserving forcing extension in which
$2^\lambda=\theta$ and $T$ has a universal model of size $\lambda^+$?
\end{Question}

\eject

\end{document}